# On alternative approach for verifiable secret sharing

*Kamil Kulesza[1], Zbigniew Kotulski[1], Joseph Pieprzyk[2]*

[1]Institute of Fundamental Technological Research, Polish Academy of Sciences, ul.Świętokrzyska 21, 00-049 Warsaw, Poland, e-mails: Kamil.Kulesza@ippt.gov.pl, Zbigniew.Kotulski@ippt.gov.pl
[2]Department of Computing, Macquarie University, NSW2109, Australia, e-mail: josef@ics.mq.edu.au

## Definitions

We propose interactive verification protocol (VP) that utilizes the concept of verification set of participants. Verification set of participants (VSoP) is the set of the secret participants that are needed to verify their shares. In order to test validity of the shares, all shares belonging to the participants from the VSoP are required.

Verification structure (VS) is the superset containing all verification sets of shares. Both terms (VSoP, VS) are closely related to authorized set of participant and general access structure , respectively. Moreover, we consider verification sets of participants that are subsets of authorized sets of participants.

For instance, consider $(t,n)$ threshold secret sharing schemes. Our proposal allows to create $(v,t,n)$ sharing scheme, where $v$ denotes number of participants in verification sets.

When participants belonging to the authorized set want to recover the secret, they first run VP for all VSoP-es contained in that set. In the example from above, VP is run for $\binom{t}{v}$ sets.

## Notation

1. Take any secret sharing scheme (SSS) over general access structure, with the $k$ participants $P_1, P_2, \ldots, P_k$ and corresponding secret shares $s_1, s_2, \ldots, s_k$. Let's denote $C_0$ as the combiner algorithm for that secret sharing scheme.

2. In order to implement VSS each secret share $s_i$ should be extended by the control part $c_i$ to form *extended secret share $s'_i$* .

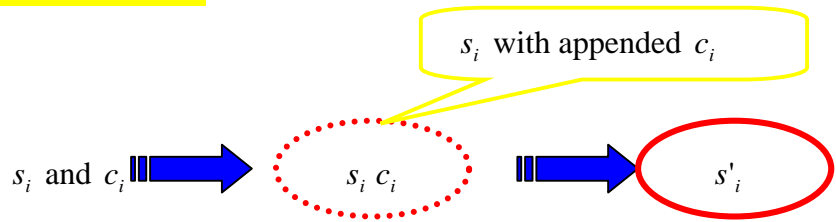

$s_i$ and $c_i$ → $s_i c_i$ → $s'_i$

($s_i$ with appended $c_i$)

3. Let $C_1(a_1,\ldots,a_\alpha), C_2(b_1,\ldots,b_\beta)$ denote the combiner algorithms for two secret shares schemes operating on the sets of the shares $a_1,\ldots,a_\alpha$ and $b_1,\ldots,b_\beta$, respectively.

## Operations in verification set of participants (VSoP)

1. Let $s'_{i_1}, \ldots, s'_{i_n}$ be extended secret shares, such that each share belongs to some participant $P_{i_\alpha}$ ($\alpha = \{1,\ldots,n\}$, $n \leq k$) and set $\{P_{i_1}, \ldots, P_{i_n}\}$ forms VSoP.

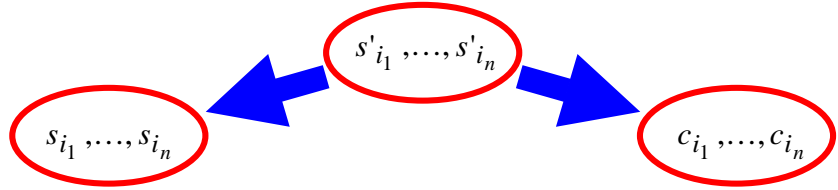

2. Combine $s_{i_1}, \ldots, s_{i_n}$ using $C_1$ to get $R_s$ (resulting secret part). Formally $C_1(s_{i_1}, \ldots, s_{i_n}) = R_s$.

3. Combine $c_{i_1}, \ldots, c_{i_n}$ using $C_2$ to get $R_c$ (resulting control part). Formally $C_2(c_{i_1}, \ldots, c_{i_n}) = R_c$

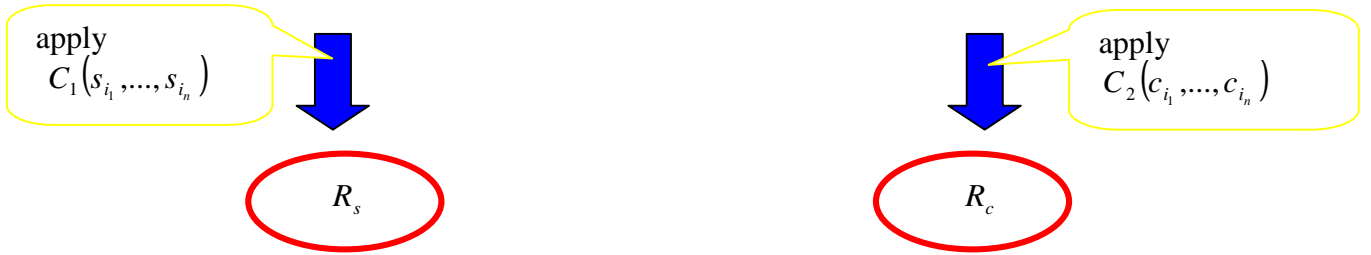

apply $C_1(s_{i_1}, \ldots, s_{i_n})$

apply $C_2(c_{i_1}, \ldots, c_{i_n})$

4. $R$ is total result equal $R_s R_c$ ($R_c$ appended to $R_s$) and set $\{P_{i_1}, \ldots, P_{i_n}\}$ forms VSoP.

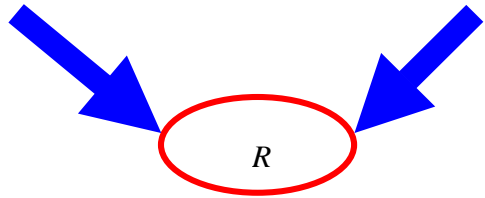

***Observation 1:*** when $s'_{i_1}, \ldots, s'_{i_n}$ are valid, then $f(R_s) = R_c$ for every $s'_i$ in the VSoP.

***Observation 2:*** different $C_1, C_2$ and $f(x)$ can be used for each VSoP.

## Construction of control function $f(x)$

*Description:*
$f(x)$ takes $l$-bit vector $x$ and computes $m$-bit image/control number, where $m < l$.

*Requirement:*
$f(x)$ should be efficient to compute.

*Sample candidates for $f(x)$:*
  a. for $m = 1$ use balanced, nonlinear Boolean function (e.g., modified bent function)
  b. for $m > 1$, one can use a vector of $m$ balanced, nonlinear Boolean functions. Consecutive values of functions from the vector are written as binary sequence to form $m$-bit control number.
  c. check-digit schemes, for instance one based on $D_5$ symmetry group.
  d. hash functions

# Verification Protocol (VP)

***Verification protocol (one round).*** The participants can verify their shares without co-operation of a third party.

1. For any verification set of shares (VSoP) compute $R$ equal $R_s R_c$.
2. Compute $f(R_s)$.
3. Test relation between $f(R_s)$ and $R_c$:

if $f(R_s) \neq R_c$ at least one of the shares in VSoP is invalid (verification is negative/ negative verification result)

if $f(R_s) = R_c$ all shares in VSoP are valid with some probability $P(s'_{i_1}, ..., s'_{i_n})$. ∎

Described above protocol is performed for all VSoP contained in the authorized set of participants, that want to recover the secret.

## On probability $P(s'_{i_1}, ..., s'_{i_n})$

Let $C_1$, $C_2$ be combiner algorithms for perfect secret sharing schemes. In addition the impact on $R$ resulting from any change of bit(s) in $s'_{i_1}, ..., s'_{i_n}$ cannot predicted in at least one of $C_1$, $C_2$.

$P(s'_{i_1}, ..., s'_{i_n})$ depends on the $R_c$ length ($m$-bits), we think that for properly chosen $f(x)$, it is related to the probability of guessing $m$-bits number. For the given $m$ there are $\left(\frac{1}{2}\right)^m$ $m$-bits numbers, hence $P(s'_{i_1}, ..., s'_{i_n}) = 1 - \left(\frac{1}{2}\right)^m$ for properly chosen $f(x)$.

## Illustrative example: $(v, t, n)$ secret sharing scheme

We assume that $f(x)$ is balanced, nonlinear Boolean function with $P(s'_{i_1}, s'_{i_2}) = \frac{1}{2}$

| $x$ | $f(x)$ |
|---|---|
| ⋮ | ⋮ |
| 00010 | 0 |
| ⋮ | ⋮ |
| 01111 | 1 |
| 10000 | 1 |
| ⋮ | ⋮ |
| 10010 | 1 |
| ⋮ | ⋮ |
| 11101 | 0 |
| ⋮ | ⋮ |
| 11111 | 0 |

Take $(3, 4)$ threshold secret sharing where secret was shared using Shamir method ($C_0$ is Shamir combiner algorithm).

Participants $P_1, P_2, P_3, P_4$ hold secret shares $s'_1, s'_2, s'_3, s'_4$ respectively.
Authorized sets of participants: $\{P_1, P_2, P_3\}, \{P_1, P_2, P_4\}, \{P_1, P_3, P_4\}, \{P_2, P_3, P_4\}$
Verification sets of participants: $\{P_1, P_2\}, \{P_1, P_3\}, \{P_1, P_4\}, \{P_2, P_3\}, \{P_2, P_4\}, \{P_3, P_4\}$
Let $g(x) = 7 + 5x + 3x^2$ be random polynomial over $GF(31)$.
$x_i = i$ for $i = 1,2,3,4$  $i \in \{1,2,3,4\}$
$s_1 = g(1) = 15 = 01111_2$, $c_1 = 1$ resulting in $s'_1 = 011111$
$s_2 = g(2) = 29 = 11101_2$, $c_2 = 0$ resulting in $s'_2 = 111010$
$s_3 = g(3) = 18 = 10010_2$, $c_3 = 1$ resulting in $s'_3 = 100101$
$s_4 = g(4) = 13 = 01101_2$, $c_4 = 1$ resulting in $s'_4 = 011011$

Let both of $C_1$, $C_2$ be combiner algorithm for KGH secret shares scheme.

Now consider authorized set $\{P_1, P_2, P_3\}$
Such authorized set has the following verification sets: $\{P_1, P_2\}, \{P_1, P_3\}, \{P_2, P_3\}$.
For:
$\{P_1, P_2\}$ $R_s = s_1 \oplus s_2 = 10010$ and $R_c = c_1 \oplus c_2 = 1$,
$\{P_1, P_3\}$ $R_s = s_1 \oplus s_3 = 11101$ and $R_c = c_1 \oplus c_3 = 0$,
$\{P_2, P_3\}$ $R_s = s_2 \oplus s_3 = 01111$ and $R_c = c_2 \oplus c_3 = 1$

*Verification protocol*

$1^{st}$ round for $\{P_1, P_2\}$ $f(R_s) = f(10010) = 1 = R_c$, hence $s'_1, s'_2$ are valid with $P(s'_1, s'_2) = \frac{1}{2}$

$2^{nd}$ round for $\{P_1, P_3\}$ $f(R_s) = f(11101) = 0 = R_c$, hence $s'_1, s'_3$ are valid with $P(s'_1, s'_3) = \frac{1}{2}$

$3^{rd}$ round for $\{P_2, P_3\}$ $f(R_s) = f(01111) = 1 = R_c$, hence $s'_2, s'_3$ are valid with $P(s'_2, s'_3) = \frac{1}{2}$

∎

*Discussion of VP results:*
1. No negative verification result was obtained in all rounds of VP.
2. Each of $s'_i$ is valid with probability $P = 1 - \left(\frac{1}{2}\right)^2 = 0{,}75$

## **Conclusions**

Presented VSS has the following features:
- works for any secret sharing scheme,
- does not require cooperation of the trusted third party,
- can be implemented for general access structure,
- its efficiency is not related to the number of dishonest participants,
- does not weaken security parameter of underlying secret sharing scheme.

The last requirement means that no extra information about the secret is revealed. For example perfect secret sharing scheme, when used with proposed VSS, still remains perfect. The information rate for the secret shares is always smaller than one, even for underlying ideal secret sharing schemes. In the particular design it can be made close to one.